\documentclass[11pt,a4paper]{article}
\usepackage{amssymb}
\usepackage{mathrsfs}
\usepackage{amsmath}
\usepackage{amstext}
\usepackage{amsfonts, amsthm, amssymb}
\usepackage{amsfonts}
\usepackage{color}
\usepackage{float}
\usepackage{graphicx}
\usepackage{subfigure}
\usepackage{caption}
\usepackage{hyperref}
\usepackage[numbers, sort&compress]{natbib}
\usepackage{amsmath}
\numberwithin{equation}{section}
\newtheorem{lem}{Lemma}[section]
\newtheorem{thm}[lem]{Theorem}

\newtheorem{claim}{Claim}[section]

\newtheorem{Proposition}[lem]{Proposition}

\begin{document}

\title{ Tight bound on the minimum degree to guarantee graphs forbidding some odd cycles to be bipartite}
\author{ Xiaoli Yuan\footnote{School of Mathematics, Hunan University,
 Changsha 410082, P. R. China. E-mail: xiaoliyuan@hnu.edu.cn.} ~~~~ Yuejian Peng\footnote{Corresponding author. School of Mathematics, Hunan University, Changsha 410082, P. R. China. E-mail: ypeng1@hnu.edu.cn. Supported in part by National Natural Science Foundation of China (No. 11931002). This paper was published on European Journal of Combinatorics 127 (2025), No. 104143. This is the final version; see https://doi.org/10.1016/j.ejc.2025.104143. }\\
}
\date{}
\maketitle
\begin{abstract}
Erd\H{o}s and Simonovits asked the following question: For an integer $r\geq 2$ and a family of non-bipartite graphs $\mathcal{H}$, determine the infimum of $\alpha$ such that any $\mathcal{H}$-free $n$-vertex graph with minimum degree at least $\alpha n$ has chromatic number at most $r$. We answer this question for $r=2$ and any family consisting of odd cycles. Let ${\mathcal C}$ be a family of odd cycles in which $C_{2\ell+1}$ is the shortest odd cycle not in ${\mathcal C}$ and $C_{2k+1}$ is the longest odd cycle in ${\mathcal C}$, we show that if $G$ is an $n$-vertex ${\mathcal C}$-free graph with $n\ge 1000k^{8}$ and $\delta(G)>\max\{ n/(2(2\ell+1)), 2n/(2k+3)\}$, then $G$ is bipartite. Moreover, the bound of the minimum degree is tight.

%Let ${\mathcal C}$ be a family of odd cycles in which $C_{2\ell+1}$ is the shortest odd cycle not in ${\mathcal C}$ and $C_{2k+1}$ is the longest odd cycle in ${\mathcal C}$.
%Let $BC_{2\ell+1}(n)$ denote the graph obtained by taking $2\ell+1$ vertex-disjoint copies of $K_{n/(2(2\ell+1)),n/(2(2\ell+1))}$ and selecting a vertex in each of them such that these vertices form a cycle of length $2\ell+1$. Note that $BC_{2\ell+1}(n)$  is ${\mathcal C}$-free  and non-bipartite with minimum degree $n/(2(2\ell+1))$. Let $C_{2k+3}(n/(2k+3))$ denote the balanced blow up of $C_{2k+3}$ with $n$ vertices.  Note that $C_{2k+3}(n/(2k+3))$ is  ${\mathcal C}$-free  and non-bipartite with minimum degree $2n/(2k+3)$.
%We show  that if $G$ is an $n$-vertex ${\mathcal C}$-free graph with $n\ge 1000k^{8}$ and  $\delta(G)>\max\{ n/(2(2\ell+1)), 2n/(2k+3)\}$, then $G$ is bipartite.  Moreover, the only $n$-vertex ${\mathcal C}$-free  non-bipartite graph with minimum degree $\max\{ n/(2(2\ell+1)), 2n/(2k+3)\}=n/(2(2\ell+1))$ is $BC_{2\ell+1}(n)$, and the only $n$-vertex ${\mathcal C}$-free  non-bipartite graph with minimum degree $\max\{ n/(2(2\ell+1)), 2n/(2k+3)\}=2n/(2k+3)$ is $C_{2k+3}(n/(2k+3))$. Our result unifies  stability results of Andr\'{a}sfai, Erd\H{o}s and S\'{o}s \cite{Andrsfi}, H\"{a}ggkvist \cite{Hggkvist} and Yuan and Peng \cite{Yuan}.
\end{abstract}
\noindent{\bf Keywords:} cycle; stability; minimum degree condition.

\noindent{\bf Mathematics Subject Classification.} 05C35.
\section{Introduction}
Let $H$ be a graph. A graph is called $H$-free if it does not contain a copy of $H$ as a subgraph. A central problem in extremal graph theory is to determine the maximum number of edges in $H$-free graphs. We use ex$(n,H)$ to denote the maximum number of edges in an $H$-free graph on $n$ vertices.
Let $T_r(n)$ denote the Tur\'{a}n graph, the complete $r$-partite graph on $n$ vertices with $r$ partition classes of size $\left\lfloor n/r\right\rfloor$ or $\left\lceil n/r\right\rceil$ and $t_r(n)=e(T_r(n))$. The classic Tur\'{a}n theorem \cite{Tur} shows that $T_r(n)$ is the unique graph attaining the maximum number of edges among all $K_{r+1}$-free graphs with $n$ vertices. Erd\H{o}s and Simonovits \cite{Erds, Simonovits} showed the stability result: if $G$ is a $K_{r+1}$-free graph with $t_r(n)-o(n^2)$ edges, then $G$ can be made into the Tur\'{a}n graph $T_r(n)$ by adding or deleting $o(n^2)$ edges in total. In place of the number of edges, it is natural to consider whether the structure of $K_{r+1}$-free graphs with large minimum degree are close to $r$-partite. In 1974, Andr\'{a}sfai, Erd\H{o}s and S\'{o}s \cite{Andrsfi} showed a seminal result as follows.
\begin{thm}[Andr\'{a}sfai, Erd\H{o}s and S\'{o}s \cite{Andrsfi}]\label{AES}
Let $r\geq 2$, if $G$ is a $K_{r+1}$-free graph with $n$ vertices and minimum degree greater than
$(3r-4)n/(3r-1)$, then $G$ is $r$-partite. Furthermore the bound is tight.
\end{thm}
%The structure of cycle in graphs with sufficiently large minimum degree has been intensively studied.
Moreover, they also proved a result on the minimum degree stability of graphs forbidding some odd cycles as follows. Let $\delta(G)$ denote the minimum degree of a graph $G$.
\begin{thm}[Andr\'{a}sfai, Erd\H{o}s and S\'{o}s \cite{Andrsfi}]\label{thmaes}
If $G$ is an $n$-vertex $\{C_3,C_5,\dots, C_{2k+1}\}$-free graph with
$$\delta(G)> \frac{2n}{2k+3},$$
then $G$ is bipartite.
\end{thm}

The existence of cycles in graphs with sufficiently large minimum degree has been intensively studied (see \cite{Bal,Brandt,Er,Hggkvist,Let,Nikiforov,Sankar,Ver}). In \cite{Hggkvist}, H\"{a}ggkvist  obtained the following related result for a single odd cycle with length no more than $9$.
\begin{thm}[H\"{a}ggkvist \cite{Hggkvist}]\label{Hag}
%Let $G$ be a non-bipartite graph with $\delta(G)>2n/(2k+3)$ and $n>\binom{k+2}{2}(2k+3)(3k+2)$. Then, if $k=\{1,2,3,4\}$, $G$ contains a $C_{2k+1}$.
Let $n>\binom{k+2}{2}(2k+3)(3k+2)$, and $G$ be an $n$-vertex $C_{2k+1}$-free graph with
$$\delta(G)>\frac{2n}{2k+3}.$$
If $k=\{1,2,3,4\}$, then $G$ is bipartite.
\end{thm}
%In 2015, F\"{u}redi and Gunderson show that the bipartite Tur\'{a}n graph attains the maximum number of edges in $C_{2k+1}$-free graphs. In 1982,
Given a graph $G$, a blow-up of a graph $G$, is obtained from $G$ by replacing
each vertex $v_i$ of $G$ with an independent set $V_i$ of  vertices
and joining each vertex in $V_i$ with each vertex in $V_j$ provided
$v_iv_j\in E(G)$. In the theorem above, the condition on the minimum degree is the best possible, this is evident by the extremal graph $C_{2k+3}(n/(2k+3))$, where $C_{2k+3}(n/(2k+3))$ denotes a blow up of $C_{2k+3}$ that each part is an independent set with size $n/(2k+3)$. In Theorem 1.3, the values of $k$ can not be extended. H\"{a}ggkvist constructed a counterexample as follows: take three vertex-disjoint copies of $K_{n/6,n/6}$, select a vertex in each of them and join every two of the selected vertices. Clearly, for $k\geq 5$, the graph is $C_{2k+1}$-free with minimum degree $n/6>2n/(2k+3)$, but it is not bipartite. Note that the minimum degree in the above construction of H\"{a}ggkvist is $n/6$.
Yuan and Peng \cite{Yuan} considered the question whether $\delta(G)>n/6$ will guarantee that an $n$-vertex $C_{2k+1}$-free graph $G$ will be bipartite for $k\ge 5$, and gave a positive answer to the question.

\begin{thm}[Yuan, Peng \cite{Yuan}]\label{thm11}
Let $k\geq 5$ and $n\ge 21000k$ be integers. Let $G$ be an $n$-vertex $C_{2k+1}$-free graph. If
$$\delta(G)\geq \frac{n}{6},$$
then either $G$ is bipartite, or $\delta(G)=n/6$ and $G$ is a graph obtained by three vertex-disjoint copies of $K_{n/6,n/6}$, selecting a vertex in each of them and joining every two of the selected vertices.
\end{thm}

What are the connections among Theorems \ref{thmaes}, \ref{Hag} and \ref{thm11}?
%Theorem \ref{thm11} implies that  $BC_{3}(n)$ (the construction of H\"{a}ggkvist)  is the unique non-bipartite $n$-vertex $C_{2k+1}$-free  graph with minimum degree $n/6$.  Note that $BC_{3}$ contains a $C_3$, we are curious what  the tight bound on the minimum degree would be to guarantee a $\{C_3, C_{2k+1}\}$-free graph to be bipartite. We discover that an $n$-vertex $\{C_3, C_{2k+1}\}$-free graph  graph $G$ with minimum degree exceeding $n/10$ will be bipartite, and the bound is the tight evident by $BC_{5}(n)$. Moreover, $BC_{5}(n)$ is the only $n$-vertex $\{C_3, C_{2k+1}\}$-free graph  graph $G$ with minimum degree  $n/10$.  Again, since $BC_{5}$ contains a $C_5$, we study what  the tight minimum degree bound would be if we add the condition that $G$ is $C_5$-free. We discover that an $n$-vertex $\{C_3, C_5, C_{2k+1}\}$-free graph  graph $G$ with minimum degree exceeding $n/14$ will be bipartite, and the bound is  tight evident by $BC_{7}(n)$. Moreover, the only $n$-vertex $\{C_3, C_5, C_{2k+1}\}$-free non-bipartite graph $G$ with minimum degree  $n/14$ is $BC_{7}(n)$.
In general, we consider what the tight bound on the minimum degree would be to guarantee an $n$-vertex  graph forbidding a family of odd cycles to be bipartite. Let ${\mathcal C}$ be a family of odd cycles. We discover that the length of the shortest odd cycle not in ${\mathcal C}$ and the length of   the longest odd cycle in ${\mathcal C}$ determine what the tight bound would be. Let $C_{2\ell+1}$ be the shortest odd cycle not in ${\mathcal C}$, and $C_{2k+1}$ be the longest odd cycle in ${\mathcal C}$.  Let $BC_{2\ell+1}(n)$ denote the graph obtained by taking $2\ell+1$ vertex-disjoint copies of $K_{n/(2(2\ell+1)),n/(2(2\ell+1))}$ and selecting a vertex in each of them such that these vertices form a cycle of length $2\ell+1$.
Observe that both  $BC_{2\ell+1}(n)$ and  the balanced blow up of $C_{2k+3}$ with $n$ vertices (denoted by $C_{2k+3}(n/(2k+3))$) are $n$-vertex ${\mathcal C}$-free  non-bipartite graphs. We discover that if $G$ is an $n$-vertex ${\mathcal C}$-free graph with $\delta(G)>\max\{ n/(2(2\ell+1)), 2n/(2k+3)\}$ (note that $\delta(BC_{2\ell+1}(n))=n/(2(2\ell+1))$ and $\delta(C_{2k+3}(n/(2k+3)))= 2n/(2k+3)$), then $G$ is bipartite. Note that graphs $BC_{2\ell+1}(n)$ and $C_{2k+3}(n/(2k+3))$ indicates that the bound is tight. Furthermore, the only $n$-vertex ${\mathcal C}$-free  non-bipartite graph with minimum degree $\max\{ n/(2(2\ell+1)), 2n/(2k+3)\}=n/(2(2\ell+1))$ is $BC_{2\ell+1}(n)$, and the only $n$-vertex ${\mathcal C}$-free  non-bipartite graph with minimum degree $\max\{ n/(2(2\ell+1)), 2n/(2k+3)\}=2n/(2k+3)$ is $C_{2k+3}(n/(2k+3))$. This is  our main result.

\begin{thm}\label{mainn}
 Let $\ell< k$ and $n\ge 1000k^{8}$ be positive integers.  Let ${\mathcal C}$ be a family of odd cycles, $C_{2\ell+1}$ be the shortest odd cycle not in ${\mathcal C}$, and $C_{2k+1}$ be the longest odd cycle in ${\mathcal C}$.
 If $G$ is an $n$-vertex ${\mathcal C}$-free graph with $\delta(G)>\max\{ n/(2(2\ell+1)), 2n/(2k+3)\}$, then $G$ is bipartite.
 Furthermore, the only $n$-vertex ${\mathcal C}$-free  non-bipartite graph with minimum degree $\max\{ n/(2(2\ell+1)), 2n/(2k+3)\}=n/(2(2\ell+1))$ is $BC_{2\ell+1}(n)$, and the only $n$-vertex ${\mathcal C}$-free  non-bipartite graph with minimum degree $\max\{ n/(2(2\ell+1)), 2n/(2k+3)\}=2n/(2k+3)$ is $C_{2k+3}(n/(2k+3))$.
 \end{thm}

 Theorem \ref{mainn} unifies Theorems \ref{thmaes}, \ref{Hag} and \ref{thm11}. Our result is also related to a question of Erd\H{o}s and Simonovits \cite{Erdos}: For an integer $r\geq 2$ and a family of non-bipartite graphs $\mathcal{H}$, what is the tight bound of $\alpha$ such that any $\mathcal{H}$-free $n$-vertex graph with minimum degree at least $\alpha n$ has chromatic number at most $r$? Theorem \ref{mainn} answers this question for $r=2$ and any family of odd cycles if $n$ is large.

 Observe that $\max\{ n/(2(2\ell+1)), 2n/(2k+3)\}=n/(2(2\ell+1))$ if $\ell<(2k-1)/8$ and $\max\{ n/(2(2\ell+1)), 2n/(2k+3)\}=2n/(2k+3)$ if $\ell>(2k-1)/8$. Also observe that a ${\mathcal C}$-free graph (${\mathcal C}$ is a family of odd cycles, where $C_{2\ell+1}$ is the shortest odd cycle not in ${\mathcal C}$, and $C_{2k+1}$ is the longest odd cycle in ${\mathcal C}$) must be $\{C_3,C_5,\dots , C_{2\ell-1};C_{2k+1}\}$-free.  Therefore, to show Theorem \ref{mainn}, it is sufficient to show the following equivalent form.

\begin{thm}\label{main}
Let $2\le \ell< k, n\ge 1000k^{8}$ be integers and let $G$ be an $n$-vertex $\{C_3,C_5,\dots , C_{2\ell-1};$ $C_{2k+1}\}$-free graph.

(i) If $ \ell>(2k-1)/8$ and $\delta(G)\ge 2n/(2k+3)$, then either $G$ is bipartite, or $\delta(G)=2n/(2k+3)$ and $G=C_{2k+3}(n/(2k+3))$.

(ii) If $2 \le \ell<(2k-1)/8$ and $\delta(G)\ge n/(2(2\ell+1))$, then either $G$ is bipartite, or $\delta(G)=n/(2(2\ell+1))$ and $G=BC_{2\ell+1}(n)$.
%$ is a graph obtained by taking $2\ell+1$ vertex-disjoint copies of $K_{n/(2(2\ell+1)),n/(2(2\ell+1))}$ and selecting a vertex in each of them such that these vertices form a cycle of length $2\ell+1$.
\end{thm}
\noindent
\textbf{Remark.} The bound $n\ge 1000k^8$ in Theorem \ref{main} can be improved by more careful calculation. Theorem \ref{Hag} and Theorem \ref{thm11} correspond to the case $\ell=1$, so we assume that $\ell\ge 2$ in Theorem \ref{main}.
%In particular, the lower bound of $n$ is not optimal, just for the convenience of calculation.

We follow standard notations through. Let $G(V,E)$ be a simple undirected graph with vertex set $V(G)$ and edge set $E(G)$, we use $e(G)$ to denote the size of $E(G)$. For $S\subseteq V(G)$, $G[S]$ denotes the subgraph of $G$ induced by $S$, we use $e(S)$ to denote the number of edges in $G[S]$. For any vertex $v\in V(G)$, $N(v)$ denote the set of all neighbors of $v$ in $G$, and $N_S(v)=N(v)\cap S$. We use $d(v)=\vert N(v)\vert$ and $d_S(v)=\vert N_S(v)\vert$. For $T\subseteq V(G)$, let $N_S(T)=\{x\in S\setminus T: \text{there exists a vertex} \ y\in T\ \text{such that} \ xy\in E(G)\}$. Let $G-T$ denote the subgraph induced by $V(G)-T$. For disjoint $X,Y\subseteq V(G)$, $G[X,Y]$ denotes the bipartite subgraph of $G$ induced by $X$ and $Y$, i.e. $G[X,Y]$ consists of all edges incident to one vertex in $X$ and one vertex in $Y$, we use
$e(X,Y)$ to denote the size of $E(G[X,Y])$. Throughout the paper, let $P_{k}$ denotes a path with $k$ vertices,  and we call an end vertex in $P_k$ an end. Let $C_{k}$ be a cycle with $k$ vertices, we call an edge that connects any two non-adjacent vertices in $C_{k}$ a chord.

\section{Proof of Theorem \ref{main}}

 In the proof, the following two results will be applied. The first one is the classical result of Erd\H{o}s and Gallai on Tur\'an numbers of paths.

\begin{thm}[Erd\H{o}s and Gallai \cite{Gallai}]\label{thm21}
For $k\geq 1$,
\begin{align*}
{\rm ex}(n,P_{k+1})\leq \frac{(k-1)n}{2}.
\end{align*}
\end{thm}
The second one guarantees that every graph $G$ contains a subgraph whose minimum degree is at least half of the average degree of $G$.
\begin{lem}\label{lem22}
If $G$ has $m$ edges and $n$ vertices, then $G$ contains a subgraph $H$ with $\delta(H)\geq \frac{m}{n}$.
\end{lem}
The above lemma has appeared in many references and can be proved by deleting vertices with minimum degree one by one till no vertex has degree smaller than half of the average degree.

Now we will give the proof of Theorem \ref{main}.

\noindent\emph{\textbf{Proof of Theorem \ref{main}}.} Let $2\le \ell< k$ and $n\ge 1000k^{8}$ be integers. Let $G$ be an $n$-vertex $\{C_3,C_5,\dots , C_{2\ell-1};C_{2k+1}\}$-free graph satisfying the minimum degree condition in Theorem \ref{main} (i) or (ii). Suppose that $G$ is not  bipartite, then there exists  an odd cycle in $G$. Let $C_{2m+1}=v_1v_2\dots  v_{2m+1}v_1$ be a shortest clockwise odd cycle of $G$. Note that $m\ge \ell\ge 2$ since $G$ is $\{C_3,\dots , C_{2\ell-1};C_{2k+1}\}$-free and $\ell\ge 2$. Let $G'=G-V(C_{2m+1})$. %By the minimality, $C_{2m+1}$ does not contain chord.
For any pair of vertices $v_i,v_j\in V(C_{2m+1})$, let $d_{C_{2m+1}}(v_i,v_j)$ denote the length of the shortest path between $v_i$ and $v_j$ in $C_{2m+1}$. Recall that $d_{C_{2m+1}}(v)$ is the number of neighbors of $v$ in $C_{2m+1}$.
\begin{claim}\label{claim21}
For any $v\in V(G')$, we have $d_{C_{2m+1}}(v)\leq 2.$ Moreover, for any $v\in V(G')$, if $v$ is adjacent to two vertices $v_i$ and $v_j$ in $C_{2m+1}$, then $d_{C_{2m+1}}(v_i,v_j)= 2$.
%$d_{C_{2m+1}}(v)\leq 3 \mbox{~for~} m=1 \mbox{~and~}$  $d_{C_{2m+1}}(v)\leq 2 \mbox{~for~} m\geq 2.$
\end{claim}
\noindent\emph{\textbf{Proof of Claim \ref{claim21}}.} %If $m=1$, then the smallest odd cycle is $C_3$. Clearly, for any vertex $v\in V(G')$, $d_{C_{3}}(v)\leq 3$.
Suppose on the contrary that there exists a vertex $x\in V(G')$ such that $d_{C_{2m+1}}(x)\geq 3$. We assume that $\{v_i,v_j,v_k\}\subseteq N_{C_{2m+1}}(x)$, where $1\leq i<j<k\leq 2m+1$. Since $m\geq 2$, any two vertices of $\{v_i,v_j,v_k\}$ are not adjacent. Otherwise, without loss of generality, we assume that $v_iv_j\in E(G)$, then $vv_iv_jv$ is a copy of $C_3$, a contradiction to $C_{2m+1}$ being a shortest cycle. Moreover, $C_{2m+1}$ is divided into three paths by $\{v_i,v_j,v_k\}$, since $C_{2m+1}$ is an odd cycle of $G$, there is at least one path whose length is odd. Without loss of generality, we assume that $v_kv_{k+1}\dots v_{2m+1}v_1\dots v_i$ is an odd path of $C_{2m+1}$. Since any two vertices of $\{v_i,v_j,v_k\}$ are not adjacent, $v_iv_{i+1}\dots v_{j}v_{j+1}\dots v_{k}$ is an path of length at least 4, %Thus, $|v_ivv_j|<|v_iv_{i-1}\dots  v_{2m+1}\dots  v_{k}\dots  v_j|$.
then we use the path $v_ivv_k$ to replace the path $v_iv_{i+1}\dots v_{j}v_{j+1}\dots v_{k}$ of $C_{2m+1}$ to get a shorter odd cycle $v_ivv_kv_{k+1}\dots v_{2m+1}v_1\dots v_i$, a contradiction. So for any vertex $v\in V(G')$, we have $d_{C_{2m+1}}(v)\leq 2.$

Suppose that $v\in V(G')$ and $v$ is adjacent to two vertices $v_i$ and $v_j$ in $C_{2m+1}$. Assume that $j>i$, if $d_{C_{2m+1}}(v_i,v_j)\neq 2$, then $d_{C_{2m+1}}(v_i,v_j)>2$, so $v_iv_{i+1}\dots v_j$ and $v_iv_{i-1}\dots v_{2m+1}\dots v_{j}$ are two paths of length at least 3. Since $C_{2m+1}$ is an odd cycle, either $v_iv_{i+1}\dots v_j$ is an odd path of $C_{2m+1}$ or $v_iv_{i-1}\dots v_{2m+1}\dots v_{j}$ is an odd path of $C_{2m+1}$. Without loss of generality, we assume that $v_iv_{i+1}\dots v_j$ is an odd path of $C_{2m+1}$, then $v_iv_{i-1}\dots v_{2m+1}\dots  v_{j}$ is an even path of $C_{2m+1}$ of length at least $4$. So $ \{v_iv_{i+1}\dots v_j\}$ together with $v$ induces a shorter odd cycle in $G$, a contradiction.
This completes the proof of Claim \ref{claim21}. $\hfill\square$
%Similarly, if $v_v_{i_1}\dots  v_{j}$ $v_jv_{j+1}\dots  v_k$

Let $S_i=N_{G'}(v_i)$ for $i \in [2m+1]$ and $S=\bigcup_{i=1}^{2m+1}S_i$. Let
$T_i=N_{G'-S}(S_i)$ for $i \in [2m+1]$,  $T=\bigcup_{i=1}^{2m+1}T_i,$ and $H=G'-S-T.$
%Now we partition $V(G')$  into three classes.

%\[
%\left\{
%\begin{array}{ll}
%S&=\bigcup\limits_{i=1}^{2m+1}S_i;\\
%T&=\bigcup\limits_{i=1}^{2m+1}T_i;\\
%H& =G'-S-T.
%\end{array}\right.
%\]
Next we show the following result giving us important structure information on  $G$.
\begin{Proposition}\label{prop23}
(i) If $ \ell>(2k-1)/8$ and $\delta(G)\ge 2n/(2k+3)$, then $m=k+1$.\\
(ii) If $2 \le \ell<(2k-1)/8$, $\delta(G)\ge n/(2(2\ell+1))$, then $m=\ell$ and for $i\neq j \in[2\ell+1]$, we have $S_i\cap S_j=\emptyset$, $T_i\cap T_j=\emptyset$,
 $E(G[S_i]) =\emptyset$, $ E(G[T_i])=\emptyset$,
 $E(G[S_i,S_j])=\emptyset, E(G[T_i,T_j])=\emptyset$, and $ E(G[S_i,T_j])=\emptyset$.
\end{Proposition}

\noindent\emph{\textbf{Proof of Proposition \ref{prop23}}.} We first show that $m$ can not be too large.
\begin{claim}\label{claim0} We have $ \ell\leq m$, moreover, the following holds. \\
(i)If $\delta(G)\geq\frac{2}{2k+3}n$, then $ m\leq k+1 (m\neq k)$.\\
(ii)If $2 \le \ell<\frac{2k-1}{8}$ and $\delta(G)\geq\frac{1}{2(2\ell+1)}n$, then $ m\leq k-1 $.
\end{claim}
\noindent\emph{\textbf{Proof of Claim \ref{claim0}}.}
 Since $G$ is $\{C_3,\dots , C_{2\ell-1};C_{2k+1}\}$-free, we have $m\ge \ell$.
Since $C_{2m+1}$ is a shortest odd cycle of $G$, it does not contain any chord. By Claim \ref{claim21}, for any $v\in V(G')$, $d_{C_{2m+1}}(v)\le 2$, so
$$(2m+1)\cdot (\delta(G)-2)\leq e(V(C_{2m+1}),V(G'))\leq 2\cdot (n-(2m+1)).$$
According to the conditions on the minimum degree in (i) and (ii), this implies that if $\delta(G)\geq 2n/(2k+3)$, then $ m\leq k+1 (m\neq k)$; and if $\ell<(2k-1)/8$ and $\delta(G)\geq n/(2(2\ell+1))$, then $ m\leq k-1 $. This completes the proof of Claim \ref{claim0}. $\hfill\square$

For $i,j\in [2m+1]$ and $i\neq j$, if we look at  the path from $v_i$ to $v_j$ (along $C_{2m+1}$ clockwisely) and the path from $v_i$ to $v_j$ (along $C_{2m+1}$ countercclockwisely), then one of them has even order (the number of vertices)  and one of them has odd order, throughout the paper, we denote the one with even order by $P_{ij}^{even}$ and its order by $p_{ij}^{even}$, and the  one with odd order by $P_{ij}^{odd}$ and its order by $p_{ij}^{odd}$.

If $m=k+1$, then Proposition \ref{prop23}(i) holds.  Assume that $m\neq k+1$, then by Claim \ref{claim0}, we have $\ell\leq m\leq k-1 $. Therefore, we only need to prove that
 under the assumption that $\ell\leq m\leq k-1 $, the assumption that $\ell>(2k-1)/8$ and $\delta(G)\geq2n/(2k+3)$ in Proposition \ref{prop23}(i) will not happen, and  
   Proposition \ref{prop23}(ii) holds.

\begin{claim}\label{claimadd}
Let $i,j\in [2m+1]$ and $i\neq j$. Then the following conclusions hold.

(i) There is no $P_{2k}$ in $G'$ with both ends in $S_i$.

(ii) There is no $P_{2k+1-p_{ij}^{odd}}$ in $G'$ with one end  in $S_i$ and another end  in $S_j$.

%(iii) There is no $P_{2k+1-p_{ij}^{even}}$ in $G'$ with one end  in $S_i$ and another end in $S_j$.

(iii) There is no $P_{2k-p_{ij}^{even}}$ in $G'$ with one end  in $S_i\setminus S_j$ and another end  in $T_j$.

(iv) There is no $P_{2k-1-p_{ij}^{odd}}$ in $G'$ with one end  in $T_i\setminus T_j$ and another end  in $T_j$.
\end{claim}
\noindent\emph{\textbf{Proof of Claim \ref{claimadd}}.} (i) If there is a $P_{2k}$ in $G'$ with both ends in $S_i$, then $V(P_{2k})$ together with $v_i$ induces a copy of $C_{2k+1}$ in $G'$, a contradiction.

(ii) If there is a $P_{2k+1-p_{ij}^{odd}}$ in $G'$ with one end  $x$ in $S_i$ and another end $y$ in $S_j$, then $V(P_{2k+1-p_{ij}^{odd}})$ together with $V(P_{ij}^{odd})$ induces a copy of $C_{2k+1}$, a contradiction.

%(iii) If there is a $P_{2k+1-p_{ij}^{even}}$ in $G'$ with one end $x$ in $S_i$ and another end $y$ in $S_j$, then this path together with edges $xv_i$, $yv_j$ and the path $P_{ij}^{even}$ forms a $C_{2k+1}$, a contradiction.

(iii) If there is a $P_{2k-p_{ij}^{even}}$ in $G'$ with one end $x$ in $S_i\setminus S_j$ and another end $y$ in $T_j$. By the definition of $T_j$, there exists a vertex $z\in S_j$ such that $yz\in E(G)$. Then $V(P_{2k-p_{ij}^{even}})\cup V(P_{ij}^{even})\cup \{z\}$ induces a copy of $C_{2k+1}$, a contradiction.

(iv) If there is a $P_{2k-1-p_{ij}^{odd}}$ in $G'$ with one end $x$ in $T_i\setminus T_j$ and another end $y$ in $T_j$. By the definition of $T_i\setminus T_j$ and $T_j$, there exist two vertices $x'\in S_i$ and $y'\in S_j$, such that $xx'\in E(G)$ and $yy'\in E(G)$. Then $V(P_{2k-1-p_{ij}^{odd}})\cup V(P_{ij}^{odd})\cup \{x',y'\}$ induces a copy of $C_{2k+1}$, a contradiction.
This completes the proof of Claim \ref{claimadd}. $\hfill\square$

The following observation is important to obtain the structure information of $G$.
\begin{claim}\label{claim22}
Let $\frac{1}{40k^5}\le \epsilon\le \frac{1}{20k^5}$. Then $|S_i\cap S_j|<\epsilon n$ for $i, j\in [2m+1]$ and $i\neq j$.
%(i)Let $\frac{1}{40k^5}\le \epsilon\le \frac{1}{20k^5}$, when $\ell>(2k-1)/8$, $\delta(G)>2n/(2k+3)$, then $|S_i\cap S_j|<\epsilon n$ for $i, j\in [2m+1]$;\\
%(ii)Let $\frac{1}{121}\le \epsilon\le \frac{1}{120}$, when $\ell<(2k-1)/8$, $\delta(G)>\frac{1}{2(2\ell+1)}n$, then $|S_i\cap S_j|<\epsilon n$ for $i, j\in [2m+1]$.
\end{claim}
\noindent\emph{\textbf{Proof of Claim \ref{claim22}}.} Suppose on the contrary that there exist $i,j\in [2m+1]$ and $i\neq j$ such that $|S_i\cap S_j|\geq\epsilon n$.
%Without loss of generality, we may assume that $j>i$,
Let $S_{ij}=S_i\cap S_j, N_T(S_{ij})=T_{ij}$, then $|S_{ij}|\geq \epsilon n$. For every vertex $v\in V(S_{ij})$, since $d(v)\geq \delta(G)$ and by Claim \ref{claim21}, we have $d_{G'}(v)\geq \delta(G)-2$.
%Let us count the number of edges incident with $S_{ij}$ in $G'$.
Then
\begin{align}\label{eq1}
|S_{ij}|\cdot (\delta(G)-2)\leq \sum\limits_{v\in S_{ij}}d_{G'}(v)\le e(S_{ij},T_{ij})+\sum\limits_{q=1}^{2m+1}e(S_{ij},S_q-S_{ij})+2e(S_{ij}).
\end{align}
By Claim \ref{claimadd}(ii), for every $q\in [2m+1]$, $G[S_{ij},S_q-S_{ij}]$ is $P_{2k+1-p_{iq}^{odd}}$-free.
By Claim \ref{claimadd}(i), $G[S_{ij}]$ is $P_{2k}$-free.  Since $m\le k-1$,
by Theorem 2.1, we have
\begin{align}\label{eq2}
\sum\limits_{q=1}^{2m+1}e(S_{ij},S_q-S_{ij})+2e(S_{ij})\leq (2m+3)kn\le (2k+1)kn.
\end{align}
Since $|S_{ij}|\geq \epsilon n$, by inequality \eqref{eq1} and \eqref{eq2}, no matter whether $\delta(G)\ge n/(2(2\ell+1))$ or $\delta(G)\ge 2n/(2k+3)$, we have
\begin{align*}
e(S_{ij},T_{ij})&\geq|S_{ij}|\cdot(\delta(G)-2)-\left(\sum\limits_{q=1}^{2m+1}e(S_{ij},S_q-S_{ij})+2e(S_{ij})\right)\\
&\geq|S_{ij}|\cdot(\delta(G)-2)-(2k+1)kn \\
&\geq \epsilon n\cdot (\delta(G)-2)-(2k+1)kn\\
&>\frac{\epsilon n}{2}\delta(G),
\end{align*}
since $\epsilon\ge \frac{1}{40k^5}$ and  $n\ge 1000k^{8}$.

By Lemma \ref{lem22}, there exists a subgraph $F\subseteq G[S_{ij},T_{ij}]$ such that $\delta(F)\geq \epsilon \delta(G)/2>2k$, thus we can greedily find a path $P_{2k+1-p_{ij}^{even}}$ in $F$ with both ends in $S_{ij}$, then $V(P_{2k+1-p_{ij}^{even}})\cup V(P_{ij}^{even})$ induces a copy of $C_{2k+1}$, a contradiction. This completes the proof of Claim \ref{claim22}.
$\hfill\square$

Let us continue the proof of Proposition \ref{prop23}. For every vertex $v_i\in V(C_{2m+1})$, $d_{G'}(v_i)\geq \delta(G)-2$, then $|S_i|\geq \delta(G)-2$ for $i\in[2m+1]$.
%And since $C_{2m+1}$ is the smallest odd cycle of $G$, $S_i$ intersects at most with two sets of $S_{i-2}$ and $S_{i+2(mod \; 2m+1)}$ in $G'$. (For convenience, we note that $S_{i+2}=S_{i+2(mod \; 2m+1)}$ as follows.)
By Claim \ref{claim22}, we have
\begin{align}\label{eq3}
|S|&\geq \sum\limits_{i=1}^{2m+1}|S_i|-\sum\limits_{\substack{i,j\in[2m+1]\\i\neq j}}|S_i\cap S_j|\notag\\
&>(2m+1)\cdot (\delta(G)-2)-\binom{2m+1}{2}\cdot \epsilon n.
\end{align}
For every $i\in[2m+1]$, let

 $$\widetilde{S_i}=S_i\setminus \bigcup_{\substack{j\in[2m+1]\setminus \{i\} }}(S_i\cap S_j),$$
 $$\widetilde{T_i}=T_i\setminus \bigcup_{\substack{j\in[2m+1]\setminus \{i\} }}(T_i\cap T_j).$$ %$|\widetilde{S_i}|=|S_i-(S_i\cap S_{i+2})-S_i\cap S_{i-2}|$ respectively.
The definitions of $\widetilde{S_i}$ and $\widetilde{T_i}$ clearly imply that $\widetilde{S_i}\cap S_q=\emptyset,\widetilde{T_i}\cap T_q=\emptyset$ for $q\in[2m+1]\backslash \{i\}$.
%$\widetilde{S_i}\cap S_l=\widetilde{S_i},\widetilde{T_i}\cap T_l=\widetilde{T_i}$ for $l=i$.
For every vertex $v\in \widetilde{S_i}$, we have $N_{C_{2m+1}}(v)=\{v_i\}$, so $d_{G'}(v)\geq \delta(G)-1$. Then

\begin{align}\label{eq4}
|\widetilde{S_i}|\cdot (\delta(G)-1)\leq \sum\limits_{v\in \widetilde{S_i}}d_{G'}(v)\le e(\widetilde{S_i},\widetilde{T_i})+\sum\limits_{\substack{q\in [2m+1]\setminus \{i\} }}e(\widetilde{S_i},S_q)+\sum\limits_{\substack{q\in [2m+1]\setminus \{i\} }}e(\widetilde{S_i},T_q)+2e(\widetilde{S_i}).
\end{align}

By Claim \ref{claimadd}(ii),  $G[\widetilde{S_i},S_q]$ is $P_{2k+1-p_{iq}^{odd}}$-free. By Claim \ref{claimadd}(iii), $G[\widetilde{S_i},T_q]$ is $P_{2k-p_{iq}^{even}}$-free. By Claim \ref{claimadd}(i),  $G[\widetilde{S_i}]$ is $P_{2k}$-free.
Thus by Theorem 2.1, we have
\begin{align}\label{eq5}
\sum\limits_{\substack{q\in [2m+1]\setminus \{i\} }}e(\widetilde{S_i},S_q)+\sum\limits_{\substack{q\in [2m+1]\setminus \{i\} }}e(\widetilde{S_i},T_q)+2e(\widetilde{S_i})\leq (4m+2)kn\le (4k-2)kn
\end{align}
since $m\le k-1$. By inequalities \eqref{eq4} and \eqref{eq5}, we have
\begin{align}\label{eq6}
e(\widetilde{S_i},\widetilde{T_i})&\geq|\widetilde{S_i}|\cdot(\delta(G)-1)-\left(\sum\limits_{\substack{q\in [2m+1]\setminus \{i\} }}e(\widetilde{S_i},S_q)+\sum\limits_{\substack{q\in [2m+1]\setminus \{i\} }}e(\widetilde{S_i},T_q)+2e(\widetilde{S_i})\right)\notag\\
&\geq |\widetilde{S_i}|\cdot (\delta(G)-1)-(4k-2)kn.
\end{align}
Let us continue the proof of Proposition \ref{prop23} according to different cases.

\textbf{Case (i) of Proposition \ref{prop23}.} $k\ge \ell>(2k-1)/8$ and $\delta(G)\geq 2n/(2k+3)$. \\
 We substitute $\delta(G)\geq 2n/(2k+3)$ into inequality \eqref{eq3} and obtain that
\begin{align}\label{eq7}
|S|
&>(2m+1)\cdot \left(\frac{2n}{2k+3}-2\right)-\binom{2m+1}{2}\cdot \epsilon n\notag\\
&=\frac{2n}{2k+3}\cdot (2m+1)-2(2m+1)-(2m+1)\cdot m\cdot \epsilon n\notag\\
&>\frac{6k}{(2k+3)\cdot 3k+1}\cdot (2m+1)n
\end{align}
since $\epsilon\le \frac{1}{20k^5}$.
By Claim \ref{claim22}, we have $|S_i\cap S_j|<\epsilon n$ for $i,j\in [2m+1]$. Moreover $|S_i|\geq 2n/(2k+3)-2$ for $i\in[2m+1]$, thus we have
\begin{equation}\label{eqsib}
|\widetilde{S_i}|=\left|S_i\setminus\bigcup_{\substack{j\in[2m+1]\setminus \{i\} }}(S_i\cap S_j)\right|>\frac{2n}{2k+3}-2-2m\cdot \epsilon n\ge \frac{n}{k+2},
\end{equation}
since $\epsilon\le \frac{1}{20k^5}$ and $n\geq 1000k^{8}$.
Substitute $\delta(G)\geq 2n/(2k+3)$ into inequality \eqref{eq6} and obtain that
\begin{align}\label{eq8}
e(\widetilde{S_i},\widetilde{T_i})&\geq |\widetilde{S_i}|\cdot \left(\frac{2n}{2k+3}-1\right)-(4k-2)kn\notag\\
&>  |\widetilde{S_i}|\cdot \frac{6kn}{(2k+3)\cdot 3k+1}
\end{align}
since $|\widetilde{S_i}|>n/(k+2)$ (see \eqref{eqsib}) and $n\ge 1000k^{8}$. Thus, we have $|\widetilde{T_i}|> 6kn/((2k+3)\cdot 3k+1)$ for every $i\in[2m+1]$. By the definition of $\widetilde{T_i}$, we have
\begin{align}\label{eq001}
|T|\geq (2m+1)\cdot|\widetilde{T_i}|>\frac{6k}{(2k+3)\cdot 3k+1}\cdot (2m+1)n.
\end{align}
By inequalities \eqref{eq7} and \eqref{eq001}, we have
$$
\frac{12k}{(2k+3)\cdot 3k+1}\cdot (2m+1)n<|S|+|T|<n.
$$
This implies that
$$
m<\frac{6k^2-3k+1}{24k}=\frac{2k-1}{8}+\frac{1}{24k}
\le \lfloor\frac{2k-1}{8}\rfloor + \frac{7}{8}+\frac{1}{24k}
< \lfloor\frac{2k-1}{8}\rfloor+1.
$$
However, it is impossible to have an integer $m$ satisfying $(2k-1)/8<m< \lfloor(2k-1)/8\rfloor+1$. Thus, when $\ell>(2k-1)/8$ and $\delta(G)\ge 2n/(2k+3)$, we have $ m\notin [\ell,k-1]$.  By Claim \ref{claim0}, we have $m=k+1$. Thus Proposition \ref{prop23}(i) holds.

Next we continue to prove Proposition \ref{prop23}(ii).

\textbf{Case (ii) of Proposition \ref{prop23}.} $\ell<(2k-1)/8$ and $\delta(G)\geq n/(2(2\ell+1))$.\\
Substituting $\delta(G)\geq n/(2(2\ell+1))$ into inequality \eqref{eq3}, we obtain that
\begin{align}\label{eq01}
|S|&\ge (2m+1)\cdot \left(\frac{n}{2(2\ell+1)}-2\right)-\binom{2m+1}{2}\cdot \epsilon n\notag\\
&=\frac{n}{2(2\ell+1)}\cdot (2m+1)-2(2m+1)-(2m+1)\cdot m\cdot \epsilon n\notag\\
&>\frac{2}{4(2\ell+1)+1}\cdot (2m+1)n
\end{align}
since $\epsilon\le \frac{1}{20k^5}$. By Claim \ref{claim22}, $|S_i\cap S_j|<\epsilon n$ for $i,j\in [2m+1]$. Moreover $|S_i|\geq n/(2(2\ell+1))-2$ for $i\in[2m+1]$, thus we have
\begin{equation}\label{eq02}
|\widetilde{S_i}|=\left|S_i\setminus\bigcup_{\substack{j\in[2m+1]\setminus \{i\}}}(S_i\cap S_j)\right|>\frac{n}{2(2\ell+1)}-2-2m\cdot \epsilon n\ge \frac{n}{2(2\ell+1)+1}
\end{equation}
since $\epsilon\le\frac{1}{20k^5}$, $m\le k-1$ and $n\ge 1000k^{8}$. By inequality \eqref{eq6}, we have
\begin{align*}
e(\widetilde{S_i},\widetilde{T_i})&\geq |\widetilde{S_i}|\cdot \left(\frac{n}{2(2\ell+1)}-1\right)-(4k-2)kn\\
&= |\widetilde{S_i}|\cdot \frac{2n}{4(2\ell+1)+1}+|\widetilde{S_i}|\cdot \left(\frac{n}{32\ell^2+36\ell+10}-1\right)-(4k-2)kn\\
&>|\widetilde{S_i}|\cdot  \frac{2n}{4(2\ell+1)+1},
\end{align*}
since $|\widetilde{S_i}|>n/(2(2\ell+1)+1)$ (see (\eqref{eq02})) and $n\ge 1000k^{8}$. Thus, we have $|\widetilde{T_i}|>2n/(4(2\ell+1)+1)$ for every $i\in[2m+1]$.  So
\begin{align}\label{eq03}
|T|\geq (2m+1)\cdot|\widetilde{T_i}|>\frac{2}{4(2\ell+1)+1}\cdot (2m+1)n.
\end{align}
By inequalities \eqref{eq01} and \eqref{eq03}, we have
$$\frac{2}{4(2\ell+1)+1}\cdot (2m+1)n\cdot 2<|S|+|T|<n.$$
This implies that $m\le \ell$. Since $G$ is $\{C_3,C_5,\dots,  C_{2\ell-1}\}$-free, we have $m=\ell$. We finish the proof for the first part of Proposition \ref{prop23}(ii) and in the rest, we will show the second part of Proposition \ref{prop23}(ii). Moreover
\begin{align}\label{eq13}
|V(H)|&=n-|S|-|T|-|V(C_{2\ell+1})|\notag\\
&<n-\frac{4n}{4(2\ell+1)+1}\cdot (2\ell+1)=\frac{n}{4(2\ell+1)+1}.
\end{align}
%By \eqref{eq8}, we have that $m=1$, i.e., the smallest odd cycle of $G$ is $C_3=\{v_1,v_2,v_3\}$.
The following claim is crucial to the proof.
\begin{claim}\label{claim23}
Let $i\in[2\ell+1]$. Then the following conclusions hold.

(i) For every vertex $v\in S_i, d_{T_i}(v)\geq 2k$;

(ii) For every vertex $v\in T_i, d_{S_i}(v)\geq 2k$.
\end{claim}
\noindent\emph{\textbf{Proof of Claim \ref{claim23}}.} (i) Suppose to the contrary that there exists a vertex $x\in S_i$ such that $d_{T_i}(x)<2k$. Then
$$d(x)-d_{T_{i}}(x)-d_{C_{2\ell+1}}(x)> \frac{n}{2(2\ell+1)}-2k-2.$$
by Pigeonhole Principle, there must exist some $S_j$ such that
\[ d_{S_j}(x)> \frac{1}{2\ell+1} \left(\frac{n}{2(2\ell+1)}-2k-2\right).\]
 Let $S_j'=N_{S_j}(x)$. Then
\begin{equation}\label{eqsj}
|S_j'|>\frac{n}{2(2\ell+1)^2}-\frac{2k+2}{2\ell+1}.
\end{equation}
Since $\delta(G)\geq n/(2(2\ell+1))$, for any vertex $v\in S_j'$, we have $d_{G'}(v)\geq n/(2(2\ell+1))-2$ by Claim \ref{claim21}. Then
\begin{equation}\label{ine0}
|S_j'|\cdot \left(\frac{n}{2(2\ell+1)}-2\right)\leq \sum\limits_{v\in S_j'}d_{G'}(v)\le e(S_j',T_j)+\sum\limits_{q=1}^{2\ell+1}e(S_j',S_q-S_j')+2e(S_j').
\end{equation}
For $q\neq j$, by Claim \ref{claimadd}(ii), $G[S_j',S_q-S_j']$ is  $P_{2k+1-p_{jq}^{odd}}$-free. For $q= j$, by Claim \ref{claimadd}(i), $G[S_j',S_j-S_j']$ is  $P_{2k}$-free. By Claim \ref{claimadd}(i),
 $G[S_j']$ is $P_{2k}$-free. Thus, by inequality (\ref{ine0}) and Theorem \ref{thm21}, we have
\begin{align*}
e(S_j',T_j)&\geq|S_j'|\cdot\left(\frac{n}{2(2\ell+1)}-2\right)-\left(\sum\limits_{q=1}^{2\ell+1}e(S_j',S_q-S_j')+2e(S_j')\right)\\
&\geq|S_j'|\cdot\left(\frac{n}{2(2\ell+1)}-2\right)-(2\ell+3)kn\\
&=|S_j'|\cdot\frac{n}{2(2\ell+1)}-(2|S_j'|+(2\ell+3)kn)\\
&=|S_j'|\cdot\frac{n}{4(2\ell+1)}+|S_j'|\cdot\frac{n}{4(2\ell+1)}-(2|S_j'|+(2\ell+3)kn)\\
&>|S_j'|\cdot \frac{n}{4(2\ell+1)},
\end{align*}
 since $|S_j'|>\frac{n}{2(2\ell+1)^2}-\frac{2k+2}{2\ell+1}$ (see (\ref{eqsj})) and $n\ge 1000k^{8}$.  Therefore,
 $$e(S_j',T_j)>|S_j'|\cdot \frac{n}{4(2\ell+1)}>\frac{n}{4(2\ell+1)^2}\cdot \frac{n}{4(2\ell+1)}=\frac{n^2}{16(2\ell+1)^3}.$$
By Lemma \ref{lem22}, there exists a subgraph $K\subseteq G[S_j',T_j]$, such that $\delta(K)>n/(16(2\ell+1)^3)>2k$. If $j=i$, then we can greedily find a path $P$ of order $2k-1$ in $K$ with both ends in $S_i'$, such that $V(P)\cup \{x,v_i\}$ induces a cycle of length $2k+1$, a contradiction. If $j\neq i$, we can greedily find a path $P_{2k-p_{ij}^{odd}}$ in $K$ with both ends in $S_j'$, such that $V(P_{2k-p_{ij}^{odd}})\cup V(P_{ij}^{odd}) \cup\{x\}$ induces a cycle of length $2k+1$, a contradiction.

%\textbf{Case 2.} $j\neq i$.

%Similarly, we obtain that
%$$e(S_j',T_j)+ \sum\limits_{l=1}^{3}e(S_j',S_l-S_j')+\sum\limits_{l=1}^{3}e(S_j',T_l-T_i)+2e(S_j')>|S_j'|\cdot (n/6-1).$$
%Then $e(S_j',T_j)>\frac{n}{36}\cdot \frac{n}{12}=\frac{n^2}{432}$. By Lemma 1, there exists a subgraph $K'\subseteq G[S_j',T_j]$, such that $\delta(K')>\frac{n}{216}$. Then we can greedily find a path $\mathcal{P}$, such that $\mathcal{P}+\{x\}$ connects an arc of $C_{3}$ to form a cycle of length $2k+1$, a contradiction.

(ii) Suppose to the contrary that there exists a vertex $y\in T_i$ such that $d_{S_i}(y)<2k$. Since $y\in T_i$, there exists a vertex $y'\in S_i$, such that $yy'\in E(G)$. Since $|H|< n/(4(2\ell+1)+1)$ (see \eqref{eq13}) and $\delta(G)\geq n/(2(2\ell+1))$, we have
$$d(y)-d_{S_i}(y)-|V(H)|>\frac{n}{2(2\ell+1)}-2k-\frac{n}{4(2\ell+1)+1}>\frac{n}{4(2\ell+1)+1},$$
since $n\geq 1000k^{8}$. According to the average principle, either there exists some $S_j$ and $j\neq i$, such that
$$d_{S_j}(y)>\frac{n}{(4(2\ell+1)+1)\cdot 2\cdot (2\ell+1)}=\frac{n}{32\ell^2+36\ell+10},$$
or there exists some $T_j$ such that $$d_{T_j}(y)>\frac{n}{32\ell^2+36\ell+10}.$$

\textbf{Case 1.} There exists some $S_j$ and $j\neq i$, such that $d_{S_j}(y)>\frac{n}{32\ell^2+36\ell+10}$.

Let $S_j''=N_{S_j}(y)$, then $|S_j''|> n/(32\ell^2+36\ell+10)$.
Since $\delta(G)\geq n/(2(2\ell+1))$, for any vertex $v\in S_j''$, we have $d_{G'}(v)\geq n/(2(2\ell+1))-2$ by Claim \ref{claim21}. Then
\begin{equation}\label{ine1}
|S_j''|\cdot \left(\frac{n}{2(2\ell+1)}-2\right)\leq \sum\limits_{v\in S_j''}d_{G'}(v)\le e(S_j'',T_j)+\sum\limits_{q=1}^{2\ell+1}e(S_j'',S_q-S_j'')+2e(S_j'').
\end{equation}
For $q\neq j$, by Claim \ref{claimadd}(ii), $G[S_j'',S_q-S_j'']$ is  $P_{2k+1-p_{jq}^{odd}}$-free. For $q= j$, by Claim \ref{claimadd}(i), $G[S_j'',S_j-S_j'']$ is  $P_{2k}$-free. By Claim \ref{claimadd}(i),
 $G[S_j'']$ is $P_{2k}$-free. Thus, by inequality (\ref{ine1}) and Theorem \ref{thm21}, we have
\begin{align*}
e(S_j'',T_j)&\geq|S_j''|\cdot\left(\frac{n}{2(2\ell+1)}-2\right)-\left(\sum\limits_{q=1}^{2\ell+1}e(S_j'',S_q-S_j'')+2e(S_j'')\right)\\
&\geq|S_j''|\cdot\left(\frac{n}{2(2\ell+1)}-2\right)-(2\ell+3)kn\\
%&\geq|S_j''|\cdotn/6-(3|S_j''|+5kn)\\
%&=|S_j''|\cdot\frac{n}{12}+|S_j''|\cdot\frac{n}{12}-(3|S_j''|+5kn)\\
&>|S_j''|\cdot \frac{n}{4(2\ell+1)},
\end{align*}
 since $|S_j''|>n/(32\ell^2+36\ell+10)$ and $n\ge 1000k^{8}$.  Therefore,
 $$e(S_j'',T_j)>|S_j''|\cdot {n \over 4(2\ell+1)}>{n\over 32\ell^2+36\ell+10}\cdot {n\over 4(2\ell+1)}={n^2\over 256\ell^3+416\ell^2+224\ell+40}.$$
By Lemma \ref{lem22}, there exists a subgraph $L\subseteq G[S_j'',T_j]$ such that $\delta(L)>n/(256\ell^3+416\ell^2+224\ell+40)>2k$ by $n\ge 1000k^8$ and $\ell\le k-1$. Then we can greedily find a path $P_{2k-1-p_{ij}^{even}}$ in $L$ with both ends in $S_j''$, such that $V(P_{2k-1-p_{ij}^{even}})\cup V(P_{ij}^{even})\cup\{y,y'\}$ induces a cycle of length $2k+1$, a contradiction.

\textbf{Case 2.} There exists some $T_j$ such that $d_{T_j}(y)>\frac{n}{32\ell^2+36\ell+10}$.

Let $T_j'=N_{T_j}(y)$. Then $|T_j'|> n/(32\ell^2+36\ell+10)$. Since $\delta(G)\geq n/(2(2\ell+1))$ and $|H|<n/(4(2\ell+1)+1)$ (see \eqref{eq13}), for any vertex $v\in T_j'$, we have
$$d_{G'-H}(v)\geq \frac{n}{2(2\ell+1)}-|H|> \frac{n}{4(2\ell+1)+1}.$$
Then
\begin{equation}\label{ine2}
|T_j'|\cdot \frac{n}{4(2\ell+1)+1}< \sum\limits_{v\in T_j'}d_{G'-H}(v)\le e(S_j,T_j')+\sum\limits_{q\neq j}e(T_j',T_q-T_j')+e(T_j',T_j-T_j')+ 2e(T_j').
\end{equation}

For $q\neq j$, by Claim \ref{claimadd}(iv), $G[T_j',T_q-T_j']$ is  $P_{2k-1-p_{jq}^{odd}}$-free. Graph $G[T_j']$ is $P_{2k}$-free since $T_j'=N_{T_j}(y)$. Then by inequality (\ref{ine2}) and Theorem \ref{thm21}, we have
\begin{align*}
e(S_j,T_j')+e(T_j',T_j-T_j')&> |T_j'|\cdot \frac{n}{4(2\ell+1)+1}-\left(\sum\limits_{q\neq j}e(T_j',T_q-T_j')+ 2e(T_j')\right)\\
&> \frac{n}{32\ell^2+36\ell+10}\cdot \frac{n}{4(2\ell+1)+1}-(2\ell+2)kn\\
&> \frac{n^2}{512\ell^3+896\ell^2+520\ell+100},
\end{align*}
since $n\geq 1000k^{8}$. According to the average principle,
either $e(S_j,T_j')> n^2/(1024\ell^3+1792\ell^2+1040\ell+200)$,
or $e(T_j',T_j-T_j')> n^2/(1024\ell^3+1792\ell^2+1040\ell+200)$.

\textbf{Case 2.1.} $e(S_j,T_j')> \frac{n^2}{1024\ell^3+1792\ell^2+1040\ell+200}$.

By Lemma \ref{lem22}, there exists a subgraph $M\subseteq G[S_j,T_j']$, such that
\[ \delta(M)>\frac{n}{1024\ell^3+1792\ell^2+1040\ell+200}>2k.\]
 If $j\neq i$, then we can greedily find a path $P_{2k-1-p_{ij}^{odd}}$ in $M$ with one end different from $y'$ in $S_j$ and another end in $T_j'$, such that $V(P_{2k-1-p_{ij}^{odd}})\cup V(P_{ij}^{odd})\cup\{y,y'\}$ induces a cycle of length $2k+1$, a contradiction. If $j=i$, then we can greedily find a path $P$ of order $2k-2$ in $M$ with one end different from $y'$ in $S_i$ and another end in $T_i'$, such that $V(P)\cup \{y,y',v_i\}$ induces a cycle of length $2k+1$, a contradiction.

\textbf{Case 2.2.} $e(T_j',T_j-T_j')> \frac{n^2}{1024\ell^3+1792\ell^2+1040\ell+200}$.

By Lemma \ref{lem22}, there exists a subgraph $G[T_j'',Q_j]\subseteq G[T_j',T_j-T_j']$, such that $\delta(G[T_j'',Q_j])>n/(1024\ell^3+1792\ell^2+1040\ell+200)$, where $T_j''\subseteq T_j'$ and  $Q_j\subseteq T_j-T_j'$.

Let $N_{S_j}(T_j'')=\{v\in S_j\setminus T_j'': \text{there exists a vertex} \ u\in T_j''\ \text{such that} \ uv\in E(G)\}.$

\textbf{Case 2.2.1} $N_{S_j}(T_j'')=\{y'\}$.

Since $\delta(G[T_j'',Q_j])>n/(1024\ell^3+1792\ell^2+1040\ell+200)>2k$, we can greedily find a path $P$ of order $2k-1$ in $G[T_j'',Q_j]$ with both ends in $T_j''$, such that $V(P)\cup\{y,y'\}$ induces a cycle of length $2k+1$, a contradiction.

\textbf{Case 2.2.2} There exists a vertex $x'\in S_j(x'\neq y')$ such that $x'\in N_{S_j}(T_j'')$.

 Since $x'\in N_{S_j}(T_j'')$, there exists a vertex $x\in T_j''$ such that $xx'\in E(G)$. If $j\neq i$. Since $\delta(G[T_j'',Q_j])>n/(1024\ell^3+1792\ell^2+1040\ell+200)>2k$, we can greedily find a path $P_{2k-2-p_{ij}^{odd}}$ in $G[T_j'',Q_j]$ with one end $x\in T_j''$ and another end in $T_j''$, such that $V(P_{2k-2-p_{ij}^{odd}})\cup V(P_{ij}^{odd})\cup \{y,y',x'\}$ induces a cycle of length $2k+1$, a contradiction. If $j=i$, we can greedily find a path $P$ of order $2k-3$ in $G[T_i'',Q_i]$ with one end $x\in T_i''$ and another end in $T_i''$, such that $V(P)\cup\{y,y',v_i,x'\}$ induces a cycle of length $2k+1$, a contradiction.
This completes the proof of Claim \ref{claim23}.   $\hfill\square$

Since we have shown that $m=\ell$ if $\ell<(2k-1)/8$ and $\delta(G)\geq n/(2(2\ell+1))$. The proof of Proposition \ref{prop23}(ii) will be completed by showing the following claim.
\begin{claim}\label{claim24}
For $i\neq j \in[2\ell+1]$, we have\\
(i)  $S_i\cap S_j=\emptyset$ and $T_i\cap T_j=\emptyset$.\\
(ii) $E(G[S_i]) =\emptyset\mbox{ and } E(G[T_i])=\emptyset$.\\
(iii) $E(G[S_i,S_j])=\emptyset, E(G[T_i,T_j])=\emptyset\mbox{ and }  E(G[S_i,T_j])=\emptyset$.
\end{claim}
\noindent\emph{\textbf{Proof of Claim \ref{claim24}}.} (i) Suppose to the contrary that there exists a vertex $x\in S_i\cap S_j$. By Claim \ref{claim23}, $G[S_i,T_i]$ has minimum degree of at least $2k$, then we can greedily find a path $P_{2k+1-p_{ij}^{even}}$ in $G[S_i,T_i]$ with one end $x$ and another end in $S_i$, such that $V(P_{2k+1-p_{ij}^{even}})\cup V(P_{ij}^{even})$ induces a cycle of length $2k+1$, a contradiction. Similarly, if there exists a vertex $y\in T_i\cap T_j$, since $y\in T_j$, there exists a vertex $y'\in S_j$, such that $yy'\in E(G)$. By Claim \ref{claim23}, we can greedily find a path $P_{2k-p_{ij}^{even}}$ in $G[S_i,T_i]$ with one end $y$ and another end in $S_i$, such that $V(P_{2k-p_{ij}^{even}})\cup V(P_{ij}^{even})\cup \{y'\}$ induces a cycle of length $2k+1$, a contradiction.
%we see that $T_i\cap T_j=\emptyset$.

(ii) Suppose to the contrary that there exists an edge $x_1y_1\in E(G[S_i])$. Then $v_ix_1y_1v_i$ forms a copy of $C_3$, a contradiction. If there exists an edge $x_2y_2\in E(G[T_i])$, since $x_2\in T_i$, there exists a vertex $x_2'\in S_i$, such that $x_2x_2'\in E(G)$. By Claim \ref{claim23}, we can greedily find a path $P$ of order $2k-2$ in $G[S_i,T_i]$ with one end $y_2\in T_i$ and another end in $S_i$, such that $V(P)\cup\{x_2,x_2',v_i\}$ induces a cycle of length $2k+1$, a contradiction.

(iii) Suppose to the contrary that there exists an edge $xy\in E(G[S_i,S_j])$, where $x\in S_i$ and $y\in S_j$. By Claim \ref{claim23}, then we can greedily find a path $P_{2k-p_{ij}^{odd}}$ in $G[S_i,T_i]$ with one end $x$ and another end in $S_i$, such that $V(P_{2k-p_{ij}^{odd}})\cup V(P_{ij}^{odd}) \cup\{y\}$ induces a cycle of length $2k+1$, a contradiction. If there exists an edge $pq\in E(G[T_i,T_j])$, where $p\in T_i$ and $q\in T_j$, then there exists a vertex $p'\in S_i$, such that $pp'\in E(G)$. By Claim \ref{claim23}, we can greedily find a path $P_{2k-1-p_{ij}^{odd}}$ in $G[S_j,T_j]$ with one end $q\in T_j$ and another end in $S_j$, such that $V(P_{2k-1-p_{ij}^{odd}})\cup V(P_{ij}^{odd})\cup\{p, p'\}$ induces a cycle of length $2k+1$, a contradiction. If there exists an edge $st\in E(G[S_i,T_j])$, where $s\in S_i$ and $t\in T_j$, then there exists a vertex $t'\in S_j$, such that $tt'\in E(G)$. By Claim \ref{claim23}, we can greedily find a path $P_{2k-1-p_{ij}^{even}}$ in $G[S_i,T_i]$ with one end $s$ and another end in $S_i$, such that
$V(P_{2k-1-p_{ij}^{even}})\cup V(P_{ij}^{even})\cup \{t,t'\}$ induces a cycle of length $2k+1$, a contradiction.
This completes the proof of Claim \ref{claim24}, consequently the proof of Proposition \ref{prop23} is complete.  $\hfill\square$

Now we will apply Proposition \ref{prop23} to prove Theorem \ref{main}. Firstly, we consider the case of Theorem \ref{main}(ii). By Proposition \ref{prop23}(ii), the shortest odd cycle of $G$ is $C_{2\ell+1}=v_1v_2\dots v_{2\ell+1}$. Since $\delta(G)\ge n/(2(2\ell+1))$, we have $|S_i|\geq n/(2(2\ell+1))-2$ and $|T_i|\geq n/(2(2\ell+1))-1$ for every $i\in[2\ell+1]$ by applying Proposition \ref{prop23}(ii). Then
\begin{align*}
|V(H)|&=n- \sum\limits_{i=1}^{2\ell+1}(|S_i|+|T_i|)-|V(C_{2\ell+1})|\\
&\leq n-(2\ell+1)\left( \frac{n}{2(2\ell+1)}-2+\frac{n}{2(2\ell+1)}-1\right)-(2\ell+1)\\
&= 2(2\ell+1).
\end{align*}

\begin{claim}\label{claim25}
$N_H(T_i)\cap N_H(T_j)=\emptyset$.
\end{claim}
\noindent\emph{\textbf{Proof of Claim \ref{claim25}}.} Suppose to the contrary that there exists a vertex $x\in N_H(T_i)\cap N_H(T_j)$. Then there exist two vertices $t_i\in T_i$ and $t_j\in T_j$, such that $xt_i,xt_j\in E(G)$, and there exists a vertex $s_j\in S_j$, such that $t_js_j\in E(G)$. By Claim \ref{claim23}, we can greedily find a path $P_{2k-2-p_{ij}^{even}}$ in $G[S_i,T_i]$ with one end $t_i$ and another end in $S_i$, such that $V(P_{2k-2-p_{ij}^{even}})\cup V(P_{ij}^{even}) \cup\{x,t_j,s_j\}$ induces a cycle of length $2k+1$, a contradiction. $\hfill\square$

\begin{claim}\label{claim26}
$|S_i|\geq \frac{n}{2(2\ell+1)}$ for every $i\in[2\ell+1]$.
\end{claim}
\noindent\emph{\textbf{Proof of Claim \ref{claim26}}.} Suppose to the contrary that there exists $S_i$ such that $|S_i|< n/(2(2\ell+1))$ for some $i\in [2m+1]$. We claim that $|S_i|+|T_i|+|N_H(T_i)|\geq n/(2\ell+1).$
For any vertex $v\in T_i$, $N(v)\subseteq S_i\cup N_H(T_i)$. Since $d(v)\geq n/(2(2\ell+1))$ and $|S_i|<n/(2(2\ell+1))$, we have
\begin{align}\label{eq9}
|N_H(T_i)|+|S_i| \geq \frac{n}{2(2\ell+1)},
\end{align}
and $N_H(T_i)\neq \emptyset$. By Claim \ref{claim25}, for any vertex $x\in N_H(T_i)$, we have $N(x)\subseteq T_i\cup V(H)$. %Since $ x\in N_H(T_i)$, there exists a vertex $t_i\in T_i$, such that $xt_i\in E(G^*)$, and there exists a vertex $s_i\in S_i$, such that $t_is_i\in E(G^*)$.%So
%\begin{align}\label{eq10}
%|T_i|\geq d(x)\geq n/6.
%\end{align}
%By \eqref{eq9} and \eqref{eq10}, we obtain
%$$|S_i|+|T_i|+|N_H(T_i)|\geq \frac{n}{3}.$$
We claim that for any vertex $x\in N_H(T_i)$, $N(x)\cap V(H)=\emptyset$. Suppose to the contrary that there exists a vertex $x'\in V(H)$ such that $xx'\in E(G)$. Since $ x\in N_H(T_i)$, there exists a vertex $t_i\in T_i$, such that $xt_i\in E(G)$, and there exists a vertex $s_i\in S_i$, such that $t_is_i\in E(G)$. Since $|H|\leq 2(2\ell+1)$, $d_T(x')\geq n/(2(2\ell+1))-2(2\ell+1)$. By Claim \ref{claim25}, there exists some $T_j$, such that $d_{T_j}(x')=d_T(x')\geq n/(2(2\ell+1))-2(2\ell+1)$. If $j=i$, then there exists a vertex $x''\in N_{T_i}(x')$. By Claim \ref{claim23}, we can greedily find a path $P$ of order $2k-4$ in $G[S_i,T_i]$ with one end $x''$ and another end in $S_i$, such that $V(P)\cup\{x',x,t_i,s_i,v_i\}$ induces a copy of $C_{2k+1}$, a contradiction. If $j\neq i$, then there exists a vertex $x'''\in N_{T_j}(x')$. By Claim \ref{claim23}, we can greedily find a path $P_{2k-3-p_{ij}^{odd}}$ in $G[S_j,T_j]$ with one end $x'''\in T_j$ and another end in $S_j$, such that $V(P_{2k-3-p_{ij}^{odd}})\cup V(P_{ij}^{odd})\cup\{x',x,t_i,s_i\}$ induces a copy of $C_{2k+1}$, a contradiction.
Thus for any vertex $x\in N_H(T_i)$, we have $N(x)\cap V(H)=\emptyset$, then $N(x)\subseteq T_i$, so
\begin{align}\label{eq10}
|T_i|\geq d(x)\geq \frac{n}{2(2\ell+1)}.
\end{align}
By inequalities \eqref{eq9} and \eqref{eq10}, we obtain that
\begin{align}\label{eq11}
|S_i|+|T_i|+|N_H(T_i)|\geq \frac{n}{2\ell+1},
\end{align}
where $|S_i|< n/(2(2\ell+1))$. If $|S_j|\geq n/(2(2\ell+1))$, then we have
\begin{align}\label{eq12}
|S_j|+|T_j|\geq \frac{n}{2(2\ell+1)}+\frac{n}{2(2\ell+1)}-1=\frac{n}{2\ell+1}-1,
\end{align}
since $|T_j|\geq n/(2(2\ell+1))-1$ for $j\in[2m+1]$. Then by inequalities \eqref{eq11} and \eqref{eq12}, we have
$V(G)=2\ell+1+\sum_{i=1}^{2\ell+1}(|S_i|+|T_i|+|N_H(T_i)|)>n$, a contradiction. This completes the proof of Claim \ref{claim26}. $\hfill\square$

By Claim \ref{claim26}, we have $|S_i|\geq n/(2(2\ell+1))$ and $|T_i|\geq n/(2(2\ell+1))-1$ for every $i\in [2\ell+1]$. Then
$$\sum\limits_{i=1}^{2\ell+1}(|S_i|+|T_i|)+|V(C_{2\ell+1})|\geq (2\ell+1)\left( \frac{n}{2(2\ell+1)}+\frac{n}{2(2\ell+1)}-1\right)+2\ell+1\geq n,$$
thus $H=\emptyset$,  $|S_i|=n/(2(2\ell+1))$ and $|T_i|=n/(2(2\ell+1))-1$ for every $i\in [2\ell+1]$, and $\delta(G)=n/(2(2\ell+1))$. In view of Proposition \ref{prop23}(ii) again, $G$ is a graph obtained by $2\ell+1$ vertex-disjoint copies of $K_{n/(2(2\ell+1)),n/(2(2\ell+1))}$, selecting a vertex in each of them such that these vertices form a cycle of length $2\ell+1$ when $\delta(G)=n/(2(2\ell+1))$.

%Now we only need to consider the case of $\delta(G)>\frac{n}{2(2\ell+1)}$. By Proposition \ref{prop23}(ii), we have for every vertex $v\in S_i$,  $d_{T_i}(v)=d(v)-d_{C_{2\ell+1}}(v)\geq \delta(G)-1$. Since $\delta(G)>\frac{n}{2(2\ell+1)}$, $\delta(G)\geq \frac{n+1}{2(2\ell+1)}$, then $|T_i|\geq\frac{n+1}{2(2\ell+1)}-1$. Moreover, by Claim \ref{claim27}, we have $|S_i|\geq \frac{n}{2(2\ell+1)}$ for $i\in [2\ell+1]$.
%Thus
%$$|S|+|T|+|V(C_{2\ell+1})|=\sum\limits_{i=1}^{2\ell+1}|S_i|+ \sum\limits_{i=1}^{2\ell+1}|T_i|+2\ell+1\geq \frac{n}{2(2\ell+1)}\cdot (2\ell+1)+\left(\frac{n+1}{2(2\ell+1)}-1\right)\cdot (2\ell+1)+2\ell+1> n.$$
%A contradiction.  So we have shown that if $\delta(G)> \frac{n}{2(2\ell+1)}$, then $G$ is bipartite.

This completes the proof of Theorem \ref{main}(ii).

Now let us prove Theorem \ref{main}(i). By Proposition \ref{prop23}(i), when $\ell>(2k-1)/8$ and $\delta(G)\geq 2n/(2k+3)$, we can conclude that $m=k+1$, i.e., the shortest odd cycle of $G$ is $C_{2k+3}=v_1v_2\dots v_{2k+3}$.
%By Claim \ref{claim0(i)}, we have the case of $m=k+1$ occurs only when $\ell>(2k-1)/8$ and $\delta(G)>\frac{2n}{2k+3}$.
By Claim \ref{claim21}, we have for any $v\in V(G')$, if $v$ is adjacent to two vertices $v_i$ and $v_j$ in $C_{2k+3}$, then $d_{C_{2k+3}}(v_i,v_j)= 2$. Let
$$D_i=N(v_{i-1})\cap N(v_{i+1}),D=\bigcup_{i\in [2k+3]}D_i.$$
By the definition of $D$ and Claim \ref{claim21}, for any two sets $D_i$ and $D_j$, we have $D_i\cap D_j=\emptyset$. So $$|D|=\sum\limits_{i\in [2k+3]}|D_i|=\sum\limits_{i,j\in [2k+3]}|N(v_i)\cap N(v_j)|.$$
According to the inclusion and exclusion principle and Claim \ref{claim21}, we have
\begin{align*}
\bigg|\bigcup_{i\in[2k+3]}N(v_i)\bigg|\geq \sum\limits_{i\in[2k+3]}|N(v_i)|-\sum\limits_{i,j\in [2k+3]}|N(v_i)\cap N(v_j)|.
\end{align*}
Since $\delta(G)\ge 2n/(2k+3)$, this implies that
$$|D|=\sum\limits_{i,j\in [2k+3]}|N(v_i)\cap N(v_j)|\geq \sum\limits_{i\in[2k+3]}|N(v_i)|-n\ge n.$$
Thus $\delta(G)=2n/(2k+3)$ and $|D|=n$, and for any $v\in V(G)$, $d_{C_{2k+3}}(v)=2$.
%If $\delta(G)>2n/(2k+3)$, we have $|D|>n$, a contradiction. Thus we have shown that when $ \ell>(2k-1)/8$ and $\delta(G)>2n/(2k+3)$, then $G$ is bipartite. %From the balanced blow up of $C_{2k+3}$, we can see that the bound $\frac{2n}{2k+3}$ is tight.
Now we are going to show that %if $G$ is a $\{C_3,C_5,\dots , C_{2\ell-1}; C_{2k+1}\}$-free non-bipartite graph on $n$ vertices with $\delta(G)=\frac{2n}{2k+3}$, then
$G$ must be a balanced blow up of $C_{2k+3}$.
%According to the above discussion, if $\delta(G)=\frac{2}{2k+3}n$, we have $|D|=n$ and for any $v\in V(G)$, $d_{C_{2k+3}}(v)=2$. We claim that
$D_i$ is an independent set for $i\in[2k+3]$. Suppose on the contrary that $D_i$ is not an independent set for some $i\in [2k+3]$, we can easily find a copy of $C_3$. Thus we conclude that $G$ is a blow up of $C_{2k+3}$. Since $\delta(G)=2n/(2k+3)$ and $|D|=n$, $|D_i|=n/(2k+3)$ for any $i\in[2k+3]$. Hence $G$ is a balanced blow up of $C_{2k+3}$. This completes the proof of Theorem \ref{main}(i).
$\hfill\square$

\section*{Remarks}

Note that Theorem \ref{thm11} implies that if $G$ is  an non-bipartite graph with minimum degree exceeding $n/6$, then $G$ contains all odd cycles with length from $11$ to $2k+1$ if $k\le n/21000$. It would be valuable to improve the requirement $n\ge 21000k$ in  Theorem \ref{thm11}. Similarly, it is interesting to improve the requirement $n\ge 1000k^{8}$ in Theorem \ref{mainn}.

In this paper, we mainly focused on studying the minimum degree stability:  For  a family of graphs $\mathcal{H}$, what is the tight bound of $\alpha$  such that the structure of any $\mathcal{H}$-free $n$-vertex graph with minimum degree at least $\alpha n$ inherits from its extremal  graph?
 Erd\H{o}s and Simonovits \cite{Erdos} asked a  question in terms of chromatic number: For an integer $r\geq 2$ and a family of  graphs $\mathcal{H}$, what is the tight bound of $\alpha$ such that any $\mathcal{H}$-free $n$-vertex graph with minimum degree at least $\alpha n$ has chromatic number at most $r$?  We do not state studies on this question in details here, but let us mention a recent interesting result of  B\"{o}ttcher, Frankl, Cecchelli, Parzcyk and Skokan \cite{Bottcher}, they proved the following: For large enough $k$, any $n$-vertex $\{C_3,C_5,\dots , C_{2k-1}\}$-free graph $G$ with minimum degree at least $n/(2k-1)$ satisfies that $\chi(G)\leq 3$. It is not clear what the tight bound on the  minimum degree condition is to guarantee the same conclusion in the above result.

\section*{Acknowledgements}
We are grateful to  both reviewers for  reading the manuscript very carefully, spending a lot of time in checking the details, and giving us valuable  comments to help improve the presentation. This work is supported by National natural science foundation of China(Nos. 11931002 and 12371327).


\begin{thebibliography}{5}
%\bibitem{allen}P. Allen, Dense $H$-free graphs are almost $(\chi(H)-1)$-partite, Electronic J.  Combinatorics, 17(1) (2010) R 21.

%\bibitem{alon} N. Alon, B. Sudakov, $H$-free graphs of large minimum degree, Electronic J. Combin. 13 (2006), R19.

\bibitem{Andrsfi}
B. Andrsf\'{a}i, P. Erd\H{o}s,  V. T. S\'{o}s,
On the connection between chromatic number, maximal clique and minimal degree of a graph,
Discrete Math. 8 (1974), 205--218.

\bibitem{Bal}
P. Balister, B. Bollob\'{a}s, O. Riordan,  R. H. Schelp,
Graphs with large maximum degree containing no odd cycles of a given length,
J. Combin. Theory Ser. B 87 (2003), 366--373.

\bibitem{Brandt}
S. Brandt, R. Faudree,  W. Goddard,
Weakly pancyclic graphs,
J. Graph Theory 27 (1998), 141--176.

\bibitem{Bottcher}
J. B\"{o}ttcher, N. Frankl, D. M. Cecchelli, O. Parzcyk, J. Skokan,
Graphs with large minimum degree and no small odd cycles are $3$-colourable, (2023),
arXiv:2302.01875v1.

\bibitem{Erds}
P. Erd\H{o}s,
On some new inequalities concerning extremal properties of graphs,
in: P. Erd\H{o}s, G. Katona (Eds.), Theory of Graphs, Academic Press, New York, 1968, pp. 77--81.

\bibitem{Er}
Erd\H{o}s, R. J. Faudree, A. Gy\'{a}rf\'{a}s, R. H. Schelp,
Odd cycles in graphs of given minimum degree, In Graph Theory, Combinatorics, and Applications,Vol. 1
(Kalamazoo, MI, 1988),  Wiley, New York, 1991, pp. 407-418.


\bibitem{Erdos}
P. Erd\H{o}s, M. Simonovits, On a valence problem in extremal graph theory, Discrete Math. 5 (1973),323--334.


\bibitem{Gallai}
P. Erd\H{o}s, T. Gallai,
On maximal paths and circuits of graphs,
Acta Math. Hungar, 10 (1959), 337--356.

\bibitem{Hggkvist}
R. H\"{a}ggkvist,
Odd cycles of specified length in nonbipartite graphs,
Graph Theory (Cambridge, 1981), North-Holland, Amsterdam, New York, 1982, pp. 89-99.

\bibitem{Kom}
J. Koml\'{o}s, M. Simonovits,
Szemer\'{e}di's regularity lemma and its applications in graph theory,
In Combinatorics, Paul Erd\H{o}s is eighty, vol. 2, J\'{a}nos Bolyai Mathematical Society, Budapest, 1996, pp. 295--352.


\bibitem{Let}
S. Letzter, R. Snyder,
The homomorphism threshold of $\{C_3, C_5\}$-free graphs,
J. Graph Theory 90 (2019), 83--106.

\bibitem{Nikiforov}
V. Nikiforov, R. H. Schelp,
Paths and cycles in graph of large minimal degree,
J Graph Theory 47 (2004), 39-52.

\bibitem{Sankar}
 M. Sankar,
 Homotopy and the Homomorphism Threshold of Odd Cycles,
 arXiv:2206.07525v1.

\bibitem{Simonovits}
M. Simonovits,
Extremal graph problems with symmetrical extremal graphs, additional chromatic conditions,
Discrete Math. 7 (1974), 349--376.

\bibitem{Tur}
P. Tur\'{a}n, Eine Extremalaufgabe aus der Graphentheorie. Mat. Fiz. Lapok 48 (1941), 436--452.

\bibitem{Ver}
J. Verstra\"{e}te, On arithmetic progressions of cycle lengths in graphs,
Combin. Probab. Comput. 9 (2000), 369--373.


\bibitem{Yuan}
X. Yuan, Y. Peng, Minimum degree stability of $C_{2k+1}$-free graphs,
J Graph Theory, (2024), 307--321.
\end{thebibliography}
\end{document}